\documentclass[a4paper,12pt]{article}
\usepackage[top=2.5cm,bottom=2.5cm,left=2.5cm,right=2.5cm]{geometry}
\usepackage{cite,amsmath,amssymb}
    \usepackage[margin=1cm,%
                font=small,%
                format=hang,%
                labelsep=period,%
                labelfont=bf]{caption}
\usepackage[algoruled,linesnumbered]{algorithm2e}
\usepackage{algorithmic}

\usepackage{amsfonts,epsf,here,graphicx}
\usepackage{latexsym,multirow,rotating}
\usepackage{pgf,tikz,subfigure}
\usepackage{epstopdf}

\newtheorem{theorem}{\bf Theorem}[section]

\newtheorem{definition}[theorem]{\bf Definition}

\begin{document}

\baselineskip=0.20in
\vspace*{10mm}

\begin{center}
{\LARGE \bf Efficient graph similarity assessment method based on vectors of topological indices}
\bigskip \bigskip

{\large \bf Matthias Dehmer$^{a,b,c,d}$,\qquad Izudin Red\v zepovi\' c$^{e}$, \\
\qquad Niko Tratnik$^{f,g}$, \qquad Petra \v Zigert Pleter\v sek$^{h,f}$
}
\bigskip\bigskip

\baselineskip=0.20in

\smallskip

$^a$ {\it Swiss Distance University of Applied Sciences, Brig, Switzerland}
\smallskip
{\tt matthias.dehmer@ffhs.ch}
\medskip

$^b$ {\it Department for Biomedical Computer Science and Mechatronics, UMIT—Private University for Health Sciences, Medical Informatics and Technology, Eduard-Wallnöfer-Zentrum 1, 6060, Hall in Tyrol, Tyrol, Austria}
\medskip

$^c$ {\it College of Artificial Intelligence, Nankai University, Tianjin, China}
\medskip

$^d$ {\it School of Science, Xian Technological University, Xian, Shaanxi 710021, China}
\medskip

$^e$ {\it State University of Novi Pazar, Serbia} \\
\smallskip
{\tt iredzepovic@np.ac.rs}
\medskip

$^f$ {\it University of Maribor, Faculty of Natural Sciences and Mathematics, Slovenia} \\
\smallskip
{\tt niko.tratnik@um.si, petra.zigert@um.si}
\medskip

$^g$ {\it Institute of Mathematics, Physics and Mechanics, Ljubljana, Slovenia}
\medskip

$^h$ {\it University of Maribor, Faculty of Chemistry and Chemical Engineering, Slovenia}


\bigskip\medskip

(Received \today)

\end{center}

\noindent
\begin{center} {\bf Abstract} \end{center}
Measuring similarity between complex objects is a fundamental task in many scientific fields. When objects are represented as graphs, graph similarity/distance measures offer a powerful framework for quantifying structural resemblance. Those comparative measures play a key role in domains such as network science, chemoinformatics, and social network analysis. While methods like graph edit distance and graph kernels are widely used, they can be computationally intensive or fail to capture fine structural variations, since they require graphs without any structural uncertainty. Another class of methods  is based on using topological indices to encode structural information of the graphs, followed by the application of  distance or similarity measures for real numbers to obtain corresponding graph-level metrics. In this paper, we introduce a novel class of distance/similarity measures which are based on multiple topological indices. Since they are generally computed in polynomial time, our method is computationally efficient in practice. We demonstrate its effectiveness through comparisons and show that it captures subtle structural information meaningfully. Additionally, we explore its applicability in two domains: analyzing random graph models in network theory and assessing molecular similarity among isomers in chemoinformatics. These preliminary results suggest that our approach holds promise for graph comparison across disciplines.
\vspace{3mm}\noindent

\baselineskip=0.30in


 \medskip\noindent
{\bf AMS Subj. Class:}  05C09, 05C82, 05C80, 05C92, 05C90, 05C12 

\section{Introduction}

{A natural outcome of observing two objects is the recognition of the degree of similarity/distance between them. 
However, accurately defining or quantifying similarity/distance is challenging due to its multidimensional nature;
it's a well-known fact that the terms similarity/distance cannot uniquely be formalized, see \cite{sobik2}.  
}

{
Given that objects from diverse domains are frequently modeled as graphs, studying graph similarity/distance measures or methods
emerges as a meaningful and relevant task. Those comparative  measures have been crucial across various fields such
as network analysis \cite{sim3}, chemoinformatics \cite{batagelj,willet_1987}, social sciences \cite{akcora} offering a 
quantitative way to assess structural resemblance between graphs \cite{sim3,wills}.
}

{
In the following, we briefly survey the most important measures in the scientific literature. The very first contribution when tackling the  
problem of measuring similarity/distance on chemical graphs has been due to Sussenguth \cite{sussenguth_1964}. 
Striktly speaking, Sussenguth dealt with exact graph matching methods, i.e., measures which are based on isomorphic relations between graphs. 
Then Zelinka \cite{zelinka1} developed the so-called Zelinka-distance; this is a isomorphism-based distance measure for graphs with the same
number of vertices which takes the size of the largest common induced subgraph of the two given graphs into account. 
Sobik \cite{sobik1} and Kaden \cite{kaden1} generalized this measure for graphs which have different number of vertices. All the mentioned approaches
so far belong to the paradigm of exact graph matching. 
}

{
A well known  graph distance metric in the field of inexact graph matching relates to the graph edit distance (GED) due to Bunke \cite{bunke1}. 
Applications thereof can be found in
\cite{gao, sanfeliu}.  Beyond graph edit distance, a variety of approaches have been developed to compare graphs more efficiently 
or to better capture structural nuances. For instance, graph kernel methods \cite{kriege} which encode structural properties via kernel functions.
}

{
Another prominent category includes measures based on topological indices or feature vectors \cite{sim0, sim1,dehmer_streib,sim2}, which 
encode structural information of the graphs by using numerical descriptors; the next step is to use existing similarity/distance measures 
for real numbers to derive comparative graph measures, see \cite{sim3}. 
%
Previous methods in this field have mostly relied on a single topological index to compute the distance or similarity between graphs. 
However, a recent approach \cite{oz} incorporates 20 different topological indices to calculate graph similarity, 
but in this  method, some information is lost when converting the vector of indices into binary representation,
which decreases its sensitivity. 
Therefore, in the present paper we aim to develop more sensitive approach that would find applications in various fields.
}

The structure of the paper is as follows. In the next section, we present our methodology. 
More precisely, distance and similarity measures based on topological indices are defined. 
Then, in Section \ref{sec3} we demonstrate that our measure, on the one hand, predicts graph similarity 
in a manner comparable to an established metric, while on the other hand, it captures subtle differences between distinct graphs more effectively. At the same time, it is significantly less computationally demanding, making it well-suited for use with larger graphs or large-scale datasets.

In the last two sections, we present preliminary results that demonstrate the potential applicability of our measure in network theory and chemistry. In Section \ref{sec4}, we analyze the behavior of our measure on several well-known random graph models, which are commonly used to simulate specific structural properties of real-world networks, such as the high clustering coefficient in the Watts–Strogatz model \cite{WS} or the presence of high-degree nodes in the Barabási-Albert preferential attachment model \cite{BA}.

Molecular similarity plays a crucial role in chemoinformatics, which has led to the development of various methods for encoding molecular structure \cite{mol_sim}. In Section \ref{sec5}, we apply our approach to assess the similarity within a set of decane isomers.


\section{{Distance and similarity measures based on topological indices}}
\label{sec2}

\subsection{Distance and similarity measures}\label{sec_sim_dist}

We start by stating some definitions regarding distance and similarity measures as they are understood in this work (see also \cite{dehmer_streib}). We do not require these measures to satisfy additional axioms (such as those of a metric). Also, we show how these measures can be converted into each other.

\begin{definition} \label{def_sim} A function  $s: X \times X \to \mathbb{R}$ is a \textbf{similarity measure} on the set $X$ if 
\begin{itemize}
\item [(1)] $0 \leq s(x,y) \leq 1$ for all $x,y \in X$,
    \item [(2)] $s(x,y) = s(y,x)$ for all $x,y \in X$, and
    \item [(3)] $s(x,x) = 1$ for all $x \in X$.
\end{itemize}
\end{definition}

\begin{definition} \label{def_dist} A function $d: X \times X \rightarrow \mathbb{R}$  is a \textbf{distance measure} (or \textbf{dissimilarity measure}) on the set $X$ if 
\begin{itemize}
\item [(1)] $d(x,y) \geq 0$ for all $x,y \in X$,
    \item [(2)] $d(x,y) = d(y,x)$ for all $x,y \in X$, and
    \item [(3)] $d(x,x) = 0$ for all $x \in X$.
\end{itemize}
\end{definition}
Additionally, if a distance measure $d$ fulfills also the following two properties:
\begin{itemize}
    \item [(4)] $d(x,y) = 0 \implies x=y$ for all $x,y \in X$, and
    \item [(5)] $d(x,z) \leq d(x,y) +  d(y,z)$ for all $x,y,z \in X$,
\end{itemize}
then we call $d$ a \textit{distance metric}. 

{
Now, let us show that we can easily get a distance measure from a similarity measure and vice versa. Firstly, we start with an arbitrary similarity measure $s$ due to Definition \ref{def_sim}. We show that a function $d: X \times X \to \mathbb{R}$, defined by \[ \forall x,y \in X,\ d(x,y) = 1- s(x,y)\] is a distance measure. 
Obviously, for any $x,y \in X$ it holds $s(x,y) \leq 1$ and therefore $d(x,y)= 1- s(x,y) \geq 0$.  
Also, we infer  \[d(x,y)=1-s(x,y)= 1-s(y,x) = d(y,x), \]
hence, the symmetry property is fulfilled.  Finally, $d(x,x) = 1- s(x,x) = 1-1 = 0$ for any $x \in X$.

This conversion can be very helpful when we start with a given similarity measure $s$ but we need to argue with distance.
Note that from $ d(x,y)= 1-s(x,y)$, we also get  $0 \leq d(x,y) \leq 1$.

For the other direction, suppose that $d$ is a distance measure on the set $X$. Firstly, suppose that there exists some constant $M > 0$ such that $d(x,y) \leq M$ for any $x,y \in X$. Then, define the function $s: X \times X \to \mathbb{R}$ as follows:
\begin{equation} \label{def_similarity}
    \forall x,y \in X,\ s(x,y) = \frac{M - d(x,y)}{M}.
\end{equation}
Choose arbitrary $x,y \in X$. Since $d(x,y) \leq M$, we get $0 \leq M - d(x,y) \leq M$ and therefore 
$$0 \leq s(x,y) = \frac{M - d(x,y)}{M} \leq 1.$$
On the other hand, 
$$s(x,y) =  \frac{M - d(x,y)}{M} =  \frac{M - d(y,x)}{M} = s(y,x).$$
Finally, we have
$$s(x,x) =  \frac{M - d(x,x)}{M} = \frac{M}{M} = 1.$$
Therefore, we have shown that $s$ is a similarity measure. 

However, if a distance measure $d$ is not bounded, then the similarity measure $s$ according to 
Definition \ref{def_dist} can be defined, for example,  by \[ \forall x,y \in X,\ s(x,y)= \frac{1}{1 + d(x,y)}. \] One can easily  check that this is a similarity measure. 

\subsection{{Graph-based distance and similarity measures}}

All graphs considered in this paper are connected. For a graph $G$, we denote by $V(G)$ the set of vertices of $G$ and by $E(G)$ the set of edges of $G$. Moreover, $d(u,v)$ is the usual shortest path \textit{distance} between vertices $u,v \in V(G)$ and $\deg(u)$ denotes the \textit{degree} of a vertex $u \in V(G)$. Additionally, the \textit{eigenvalues} of a graph $G$ are the eigenvalues of the adjacency matrix of $G$.

Now, we present a new method for computing the similarity between two graphs based on topological indices.   
 Let $I_1, \ldots, I_k$ be topological indices, where $k \geq 3$. For a graph $G$, let $M(G) = \max \{ I_1(G), \ldots, I_k(G) \}$  and $m(G) = \min \{ I_1(G), \ldots, I_k(G) \}$. Since topological indices can have very different values, we begin by scaling them to enable meaningful comparison. For any $i \in \{ 1, \ldots, k \}$, define
 $$J_i (G) = \frac{ I_i(G) - m}{M-m}. $$
 Note that $J_i(G) \in [0,1]$ for any $i \in \{ 1, \ldots, k \}$. Then, for a given $p \geq 1$ one can calculate the distance between two graphs, $G_1$ and $G_2$, as
$$d_p(G_1,G_2) = d_p \Big( \big( J_1(G_1), \ldots, J_k(G_1) \big),    \big(J_1(G_2), \ldots, J_k(G_2) \big) \Big),$$
where $d_p$ is a $p$-distance on $\mathbb{R}^k$, defined as
$$d_p \big( (x_1, \ldots, x_k), (y_1, \ldots, y_k) \big) = \left( \sum_{i=1}^k |x_i - y_i|^p \right)^{\frac{1}{p}}$$
for any two $(x_1, \ldots, x_k), (y_1, \ldots, y_k) \in \mathbb{R}^k$.

Note that $0 \leq d_p(G_1,G_2) \leq k^{\frac{1}{p}}$ since $J_i(G_\ell) \in [0,1] $ for  $i \in \{ 1, \ldots, k \}$ and $\ell \in \{ 1,2 \}$. If 
$p=1$, one obtains the \textit{Manhattan distance}:
 $$d_1(G_1,G_2) = \sum_{i=1}^k \big| J_k(G_1) - J_k(G_2) \big|.$$

Moreover, if $p=2$, we get the \textit{Euclidean distance}: 
 \begin{equation*} 
 d_2(G_1,G_2) = \sqrt{ \sum_{i=1}^k (J_i(G_1) - J_i(G_2))^2}. 
 \end{equation*}

We chose the $d_p$ distance because it is a standard, well-understood measure for comparing vectors, offering both computational efficiency and flexibility through the choice of $p$. Note that according to Definition \ref{def_dist}, $d_p$ is obviously a distance measure for any chosen set of graphs. 


Finally, we define the similarity between two graphs $G_1$ and $G_2$ by employing formula \eqref{def_similarity} as follows:
$$s_p(G_1,G_2) = \frac{k^{\frac{1}{p}} - d_p(G_1,G_2)}{k^{\frac{1}{p}}}.$$
{
Based on the discussion in Section \ref{sec_sim_dist}, it is obvious that $s_p$ is a similarity measure.
}


We observe that for a fixed $k$, the time complexity of the proposed similarity measure is asymptotically equal to the time complexity of the most computationally expensive topological index used in its computation.

Of course, the obtained similarity measure depends on the chosen topological indices $I_1, \ldots, I_k$. Therefore, in the selection of topological indices we decided to include one distance-based index, one degree-based index, two distance-degree-based indices, and two eigenvalue-based indices. Let $G$ be a connected graph on $n$ vertices with eigenvalues $\lambda_1, \ldots, \lambda_n$. The following indices were used:
\begin{enumerate}
    \item the \textit{Harary index} \cite{ivanciuc,plavsic}, defined as $$\displaystyle H(G) = \sum_{ \substack{ \{ u, v \} \subseteq V(G) \\ u \neq v}} \frac{1}{d(u,v)},$$
    \item the \textit{Sombor index} \cite{sombor}, defined as $$\displaystyle SO(G) = \sum_{uv \in E(G)} \sqrt{(\deg (u))^2 + (\deg (v))^2},$$
    \item the \textit{degree distance} \cite{dobrynin} (also known as the \textit{Schultz index} \cite{schultz}), defined as $$\displaystyle DD(G) = \sum_{ \{ u , v \} \subseteq V(G)} (\deg (u) + \deg (v) ) d(u,v),$$
    \item the \textit{Gutman index} \cite{gutman_index}, defined as $$\displaystyle Gut(G) = \sum_{ \{ u , v \} \subseteq V(G)} \deg (u)  \deg (v)  d(u,v),$$
    \item the \textit{graph energy} \cite{energy}, defined as 
    $$E(G) = \sum_{i=1}^n |\lambda_i|,$$
    \item the \textit{Estrada index} \cite{estrada}, defined as 
    $$EE(G) = \sum_{i=1}^n e^{\lambda_i}.$$
\end{enumerate}

Indices were chosen based on their prevalence in the literature, with degeneracy \cite{konst}  taken into consideration. Note that degeneracy of a topological index measures its power to distinguish between graphs. For some  results on this topic see \cite{root-wiener,DGF12,izudin,izudin_furtula}. 

One of the main characteristics of any computational technique is time complexity, so we analyze this property of our method. Let $G_1$ and $G_2$ be connected graphs with $n_1$ and $n_2$ vertices, respectively.  Denote $n = \max \{ n_1, n_2 \}$. The distances between all pairs of vertices in a connected graph on $n$ vertices can be computed in ${\cal{O}}\left(n^3\right)$ time, for example by using the Floyd–Warshall algorithm. Additionally, the degrees of all vertices can be computed in ${\cal{O}}\left(n^2\right)$ time when the graph is represented by the adjacency matrix. Consequently, the Harary index, the degree-distance, and the Gutman index can be calculated in ${\cal{O}}\left(n^3\right)$ time, while  the Sombor index takes ${\cal{O}}\left(n^2\right)$ time.

On the other hand, approximations of the eigenvalues of a graph can be computed, for example, by using the QR algorithm \cite{francis}, where each iteration requires ${\cal{O}}\left(n^3\right)$ time (however, since the adjacency matrix is symmetric, the time complexity of an individual step could be reduced  after transforming the matrix into another form \cite{demmel,ortega}). Therefore, approximating the eigenvalues using the standard QR algorithm takes at most ${\cal{O}}\left(sn^3\right)$ time, where $s$ denotes the number of iterations needed to achieve the desired numerical precision. Consequently, the same upper bound on time complexity applies to the computation of the graph energy, the Estrada index, and the similarity measure $s_p(G_1,G_2)$, where $p \geq 1$ and  $I_1, \ldots, I_6$ are the topological indices mentioned above.

In order to test the dependence of our approach on the number of indices used, 
we expanded the selection of the above mentioned indices with the first Zagreb index \cite{gutman_zagreb}, the Randi\' c index \cite{randic}, the resolvent energy \cite{resolvent}, and the Wiener index \cite{wiener}. Therefore,  the number of components of vectors in this approach equals $k=10$. Figure \ref{corr} shows the correlation between the similarity $s_2$ calculated by using 6 indices and the similarity $s_2$ calculated by using 10 indices for the family of all trees with 7 vertices.  We can see that the correlation is very strong ($R= 0.9992$), so the larger number of indices does not significantly affects the similarity results. Therefore, we continue the investigation with only 6 indices. 

\begin{figure}[h!] \centering
    \includegraphics[scale=.7]{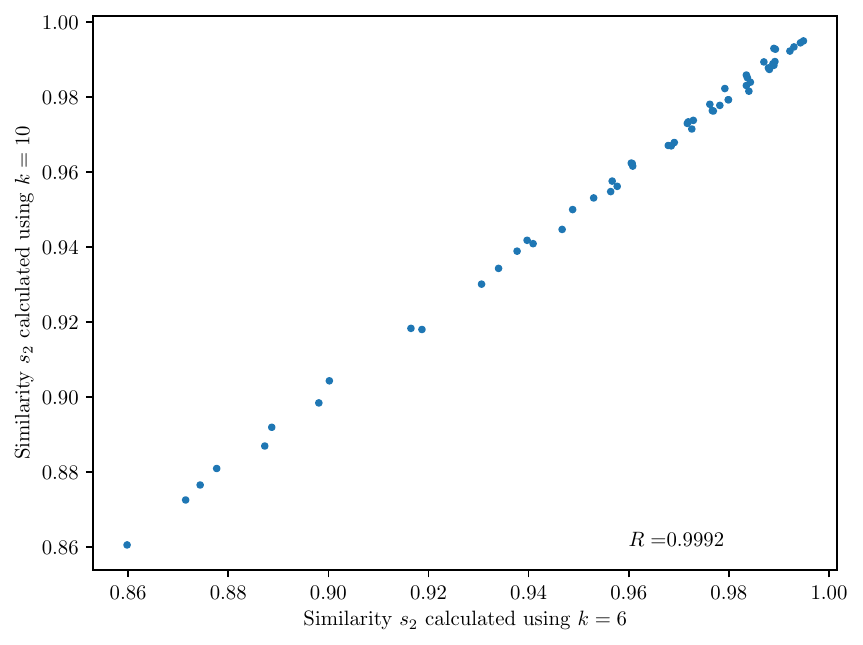}  
     \caption{Correlation between similarity $s_2$ results obtained by $k=6$ and $k=10$ for the family of all trees with 7 vertices.}
     \label{corr}
 \end{figure}

\section{Numerical results}
\label{sec3}

In this section, we present several computational results related to graph similarity measures $s_1$ and $s_2$ defined in the previous section. We consider the family $T_{7}$ of all trees on 7 vertices and the family $N_7$ of all connected graphs on 7 vertices. Note that $T_{7}$ contains 11 trees and $N_7$ contains 853 graphs. Moreover, if a graph family contains $m$ graphs, then there are $\binom{m}{2}$ pairs to consider. For the family $T_{7}$ that means 55 pairs and for the family $N_7$ that means 363378 pairs of graphs.

We compare similarity measures $s_1$ and $s_2$ to a similarity measure obtained by the \textit{graph edit distance} ($GED$). For graphs $G_1$ and $G_2$, ${GED}(G_1,G_2)$ is defined as the minimum number of elementary graph operations that are needed to transform graph $G_1$ into  graph $G_2$ \cite{sanfeliu}. This is a widely used measure for quantifying the similarity between two graphs, but computing it is generally NP-hard, making it computationally expensive for large graphs. In  \cite{gao}, a review on the graph edit distance can be found. Note the similarity $s_{GED}$ between graphs $G_1$ and $G_2$ can then be defined, for example, as
$$s_{GED} = \frac{1}{GED(G_1,G_2) + 1}.$$

Firstly, we present two figures on which one can see the similarity results for all  pairs of graphs, which are ordered according to their similarity, see Figure \ref{all_pairs}. As already mentioned, the computation of graph edit distance is highly costly, therefore we managed to compute it only for $T_7$. 

\begin{figure}[h!] 
\centering
    \includegraphics[width=.5\textwidth]{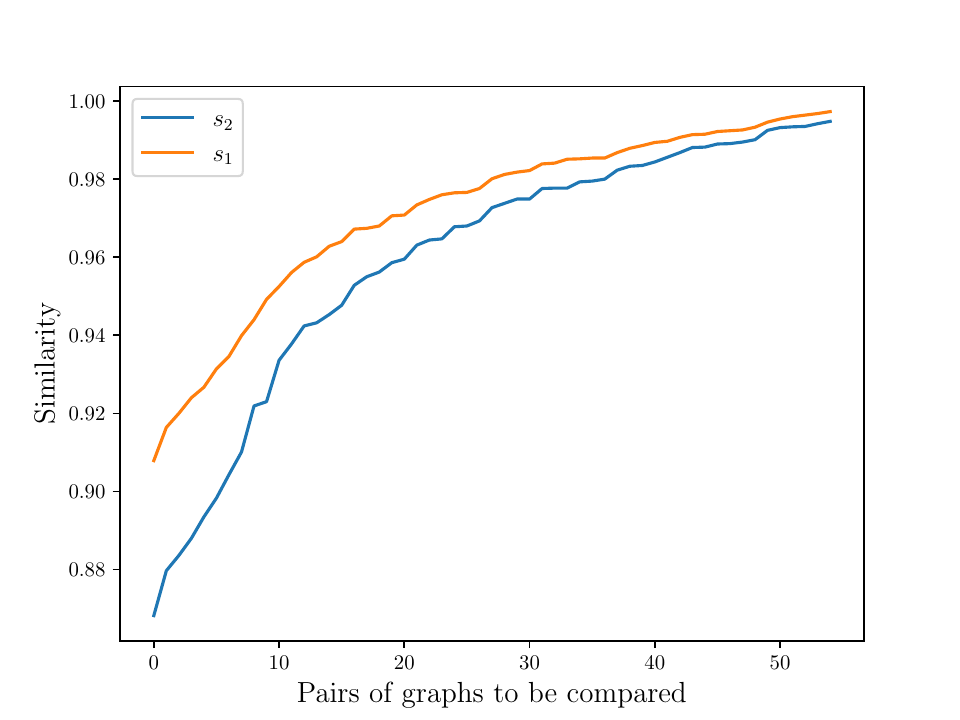} \hspace{-2em}
    \includegraphics[width=.5\textwidth]{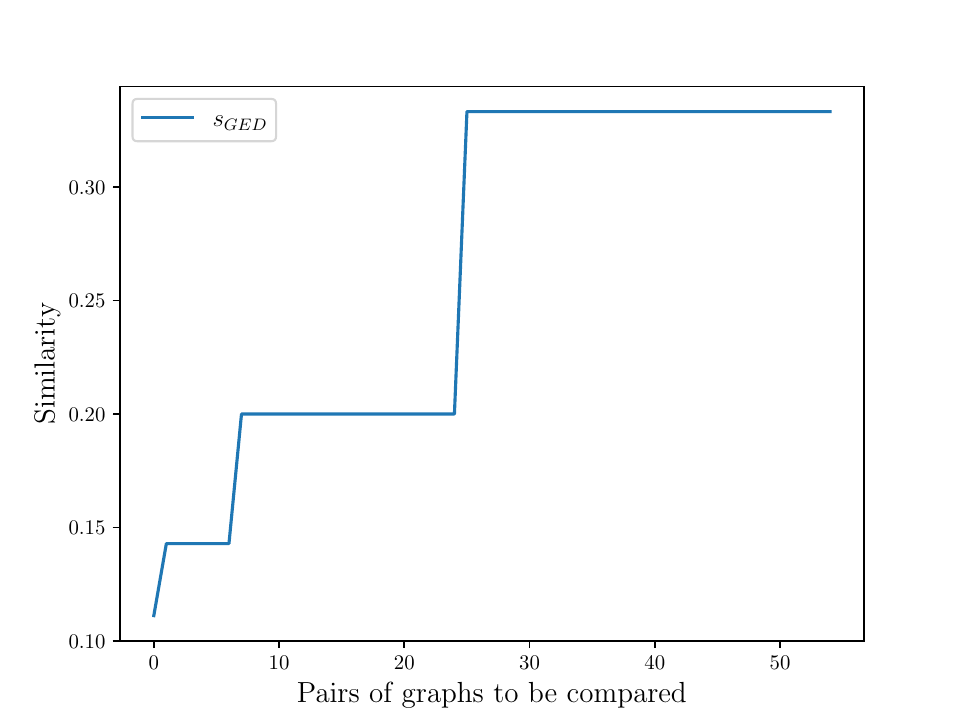} 
     \includegraphics[width=.5\textwidth]{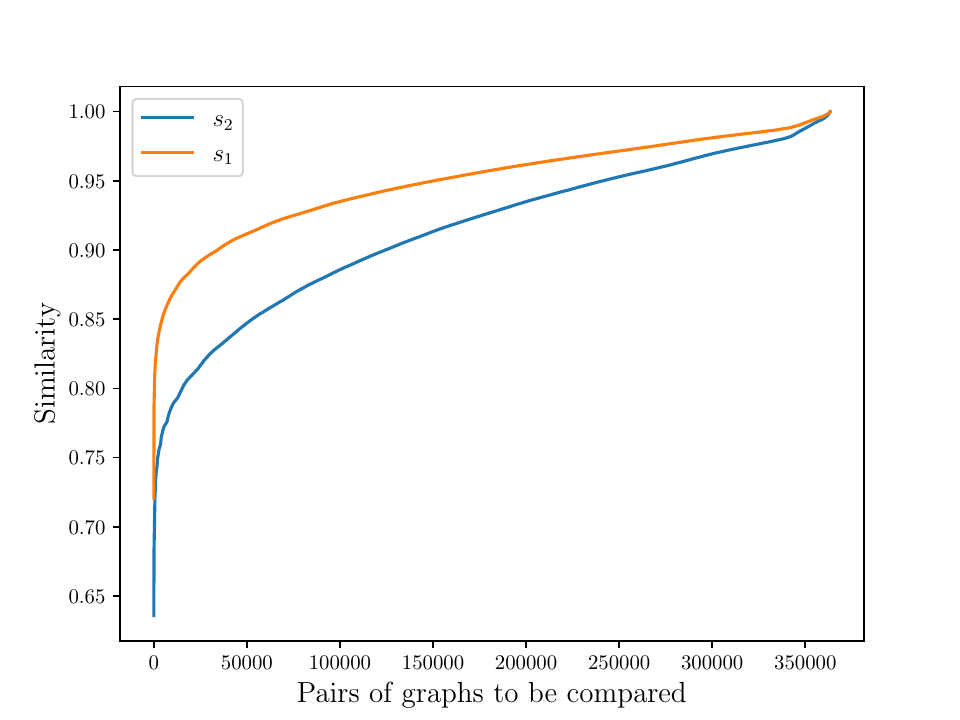}
     \caption{Similarity for all pairs of graphs in families $T_{7}$ (first row) and $N_7$ (second row).}
     \label{all_pairs}
 \end{figure}

One can observe different behavior of our measures compared to the similarity based on graph edit distance. In particular, $s_1$ and $s_2$
  distinguish more effectively between pairs of graphs, as shown by the smoothly increasing trends in Figure \ref{all_pairs}. 

Moreover, we examine correlations between similarity measures $s_1$, $s_2$, and $s_{GED}$, see Figure \ref{corr12}. We notice that the results obtained by similarities $s_1$ and $s_2$ are comparable. Both of them better discriminate between pairs of graphs in comparison to $s_{GED}$. However, they show the same trend since with increasing $s_{GED}$ also $s_1$ and $s_2$ increase. Based on these findings, we see that it is enough to proceed with $s_2$. 

\begin{figure}[h!]
\centering
    \hspace{-2em} \includegraphics[width=.48\textwidth]{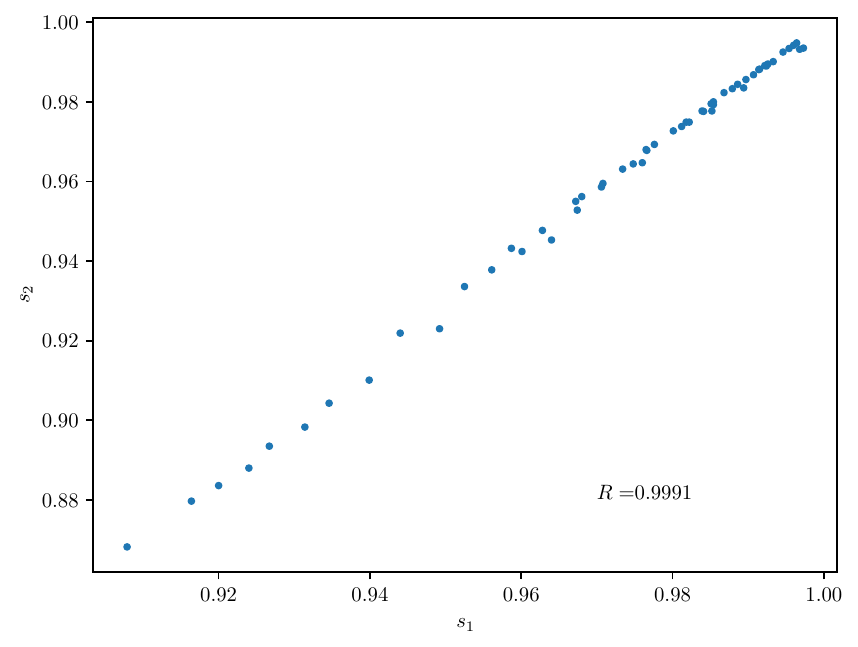} \\
    \includegraphics[width=.5\textwidth]{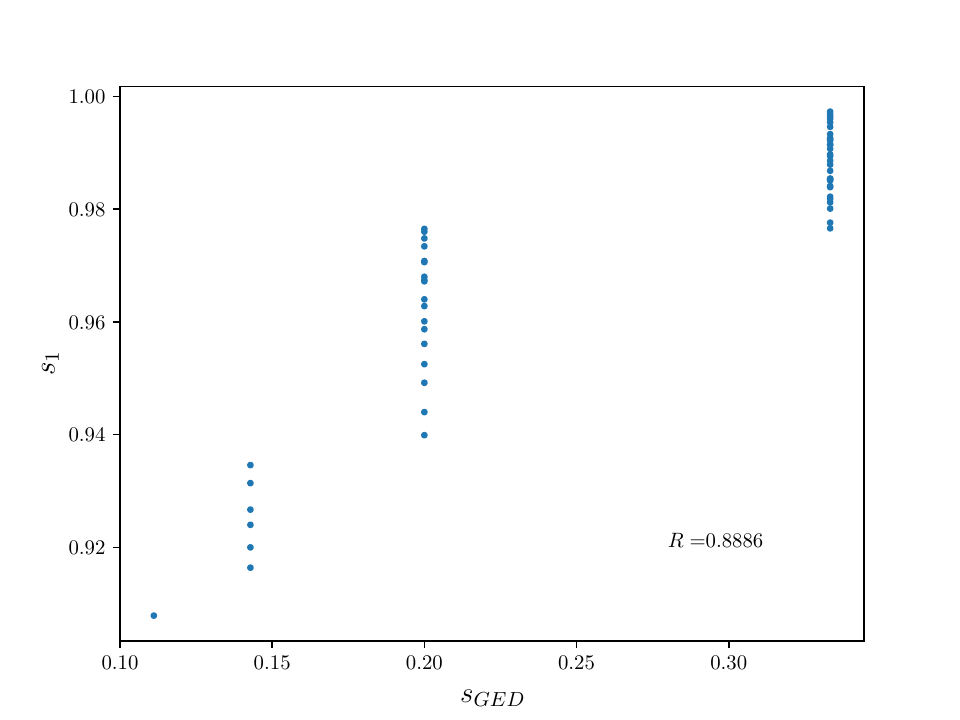} \hspace{-2em}
    \includegraphics[width=.5\textwidth]{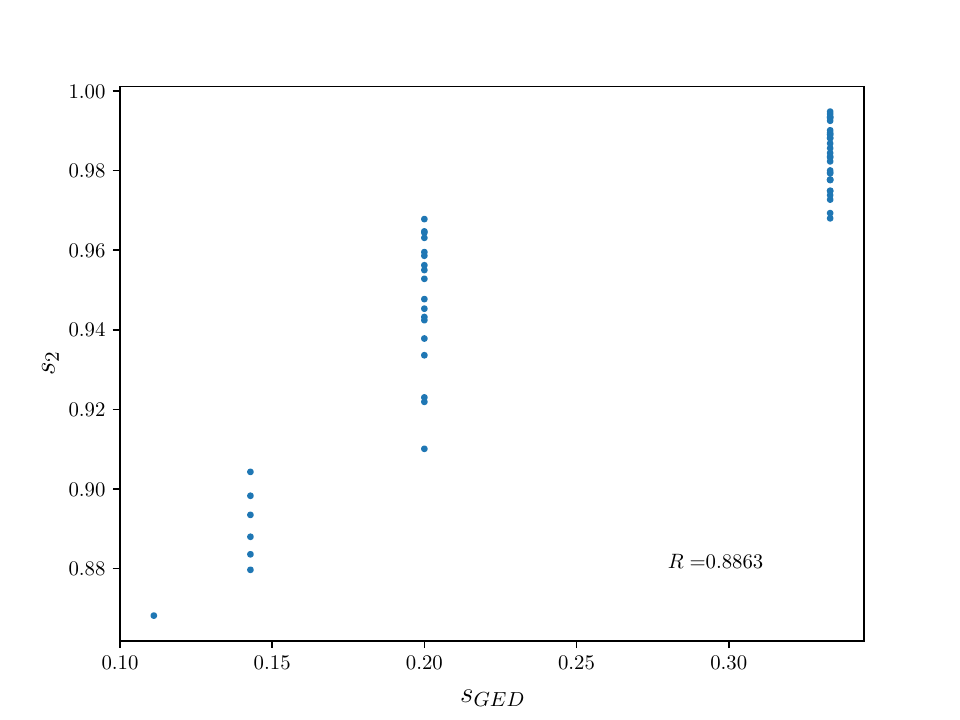}   
     \caption{Correlations between $s_1$, $s_2$ and $s_{GED}$ in case of $T_{7}$.}
     \label{corr12}
 \end{figure}

Finally, to evaluate the similarity results, we examine the pairs of graphs with minimum and maximum similarity for both families of graphs. Graphs with the minimum and maximum similarity $s_2$ in family $T_{7}$ are depicted in Figure \ref{min_trees7} and Figure \ref{max_trees7}, respectively. Note that there is only one pair of trees that achieve minimum, while on the other hand there several pairs with the maximum similarity. Such behavior occurs also in the family $N_7$. Graphs with the minimum and maximum similarity $s_2$ in this family are depicted in Figure \ref{min_graphs7} and Figure \ref{max_graphs7}, respectively.

\begin{figure}[h!]
\centering
     \includegraphics[width=.3\textwidth]{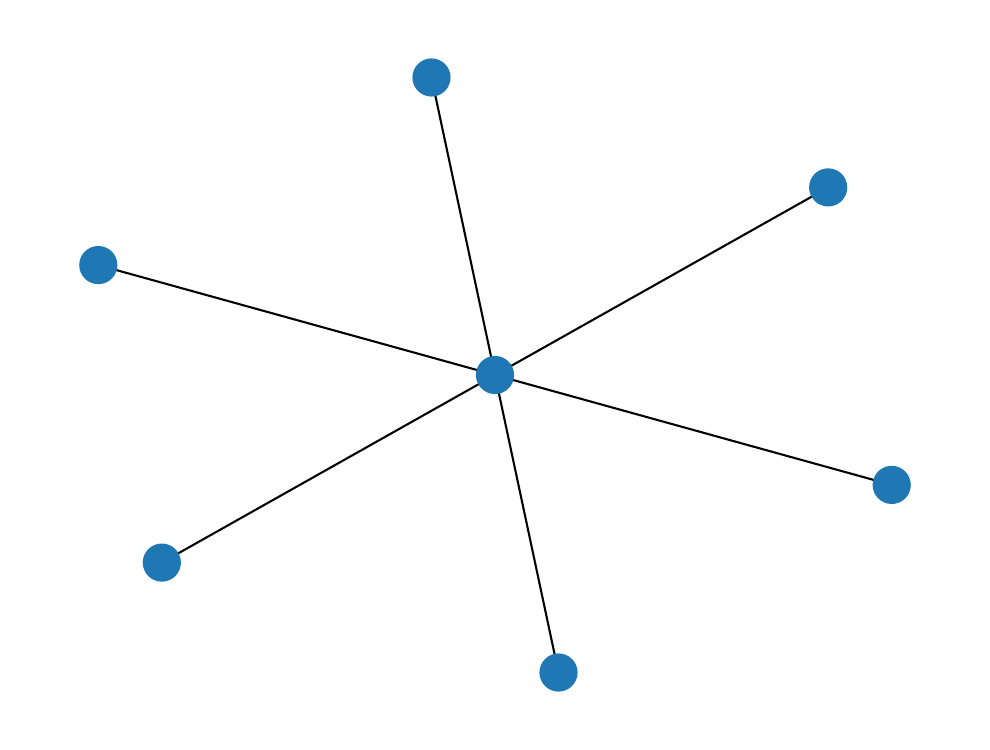} \ \ \ \ \ \
    \includegraphics[width=.3\textwidth]{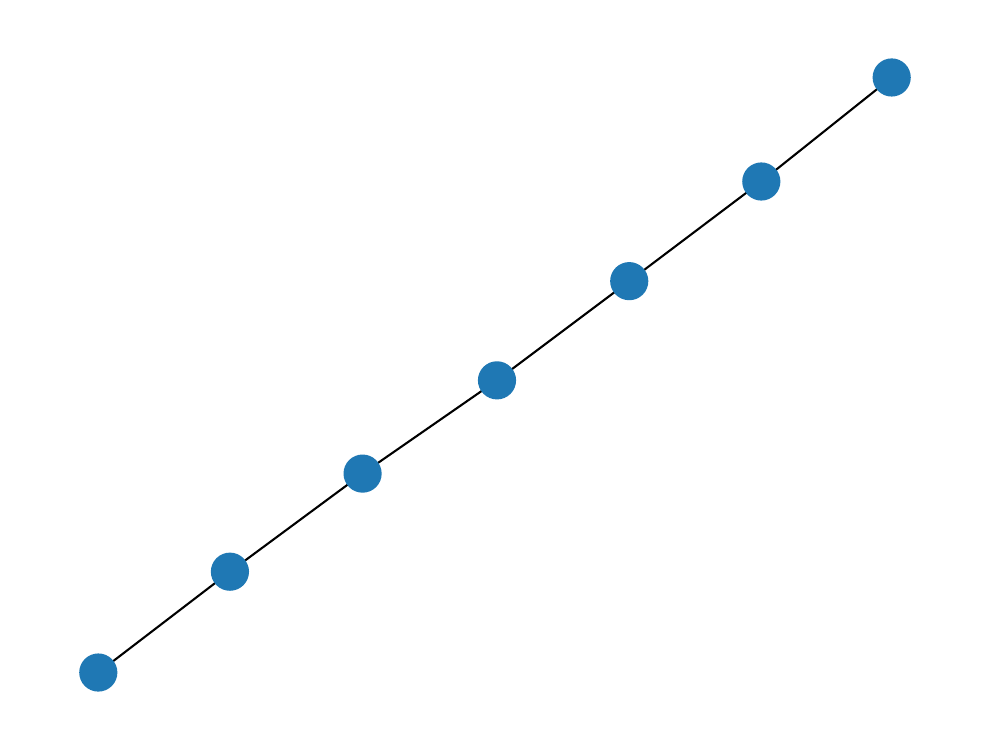}   
     \caption{Uniquely determined pair of graphs with the minimum similarity $s_2$ in family $T_{7}$.}
     \label{min_trees7}
 \end{figure}

 \begin{figure}[h!]
\centering
     \includegraphics[width=.3\textwidth]{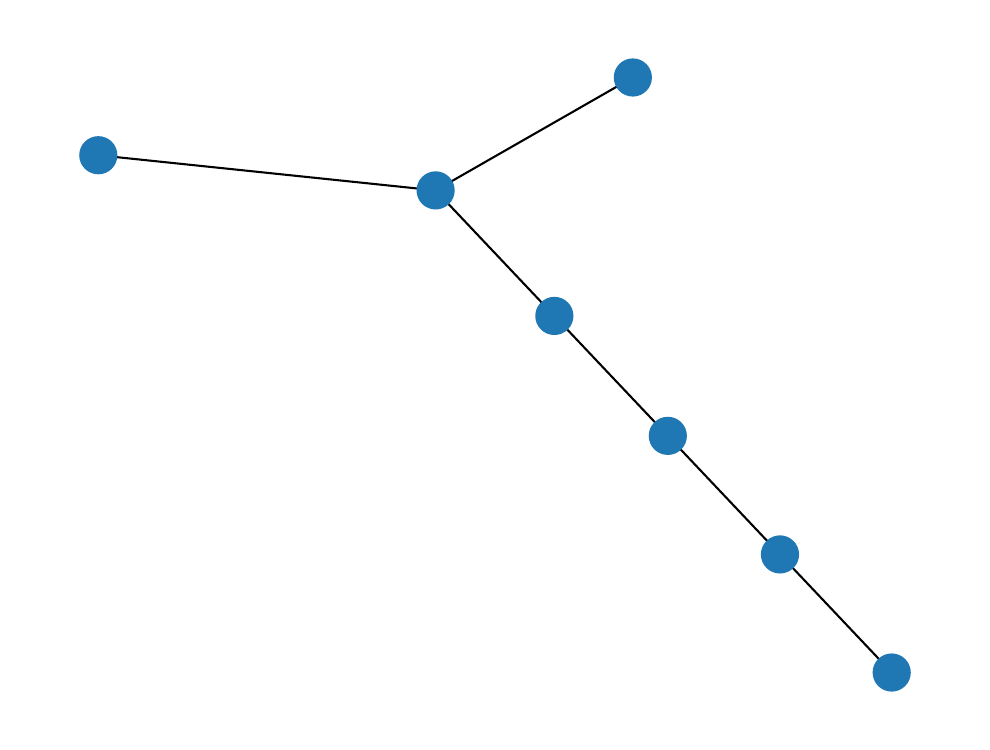} \ \ \ \ \ \
    \includegraphics[width=.3\textwidth]{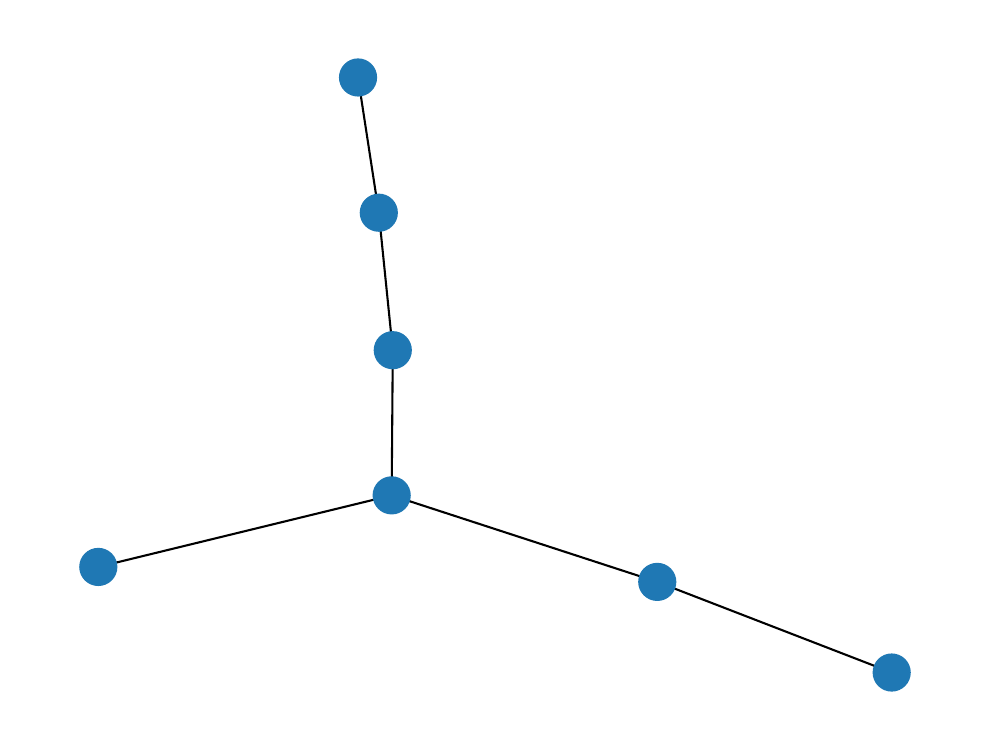}   
     \caption{One pair of graphs with the maximum similarity $s_2$ in family $T_{7}$.}
     \label{max_trees7}
 \end{figure}

 \begin{figure}[h!]
\centering
     \includegraphics[width=.3\textwidth]{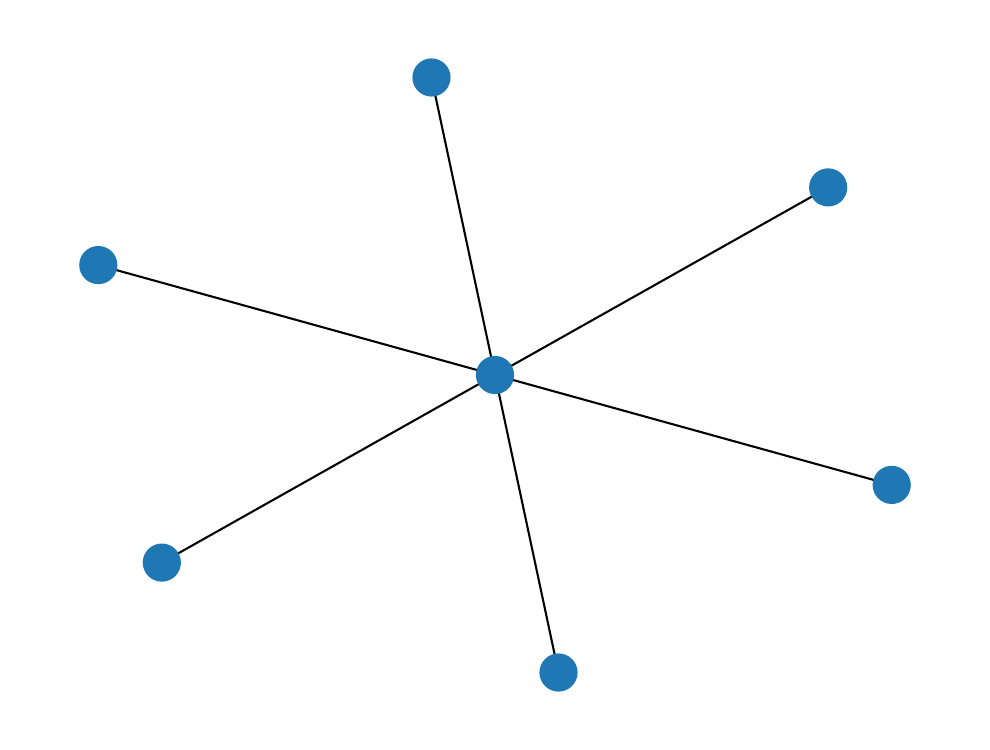} \ \ \ \ \ \
    \includegraphics[width=.3\textwidth]{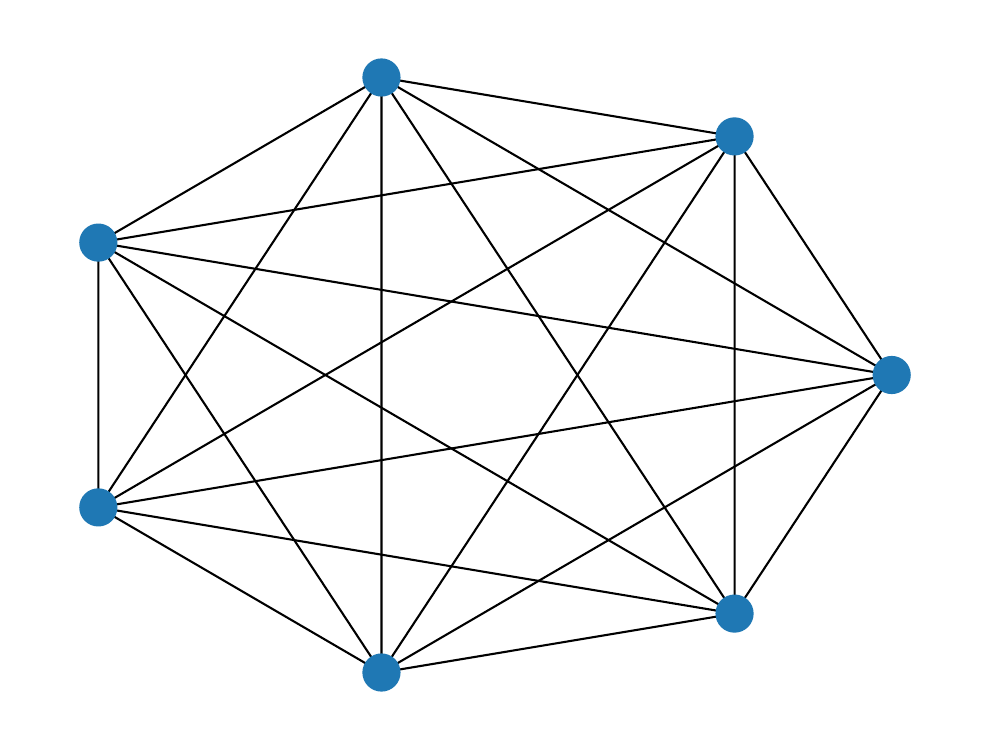}   
     \caption{Uniquely determined pair of graphs with the minimum similarity $s_2$ in family $N_7$.}
     \label{min_graphs7}
 \end{figure}

 \begin{figure}[h!]
\centering
     \includegraphics[width=.3\textwidth]{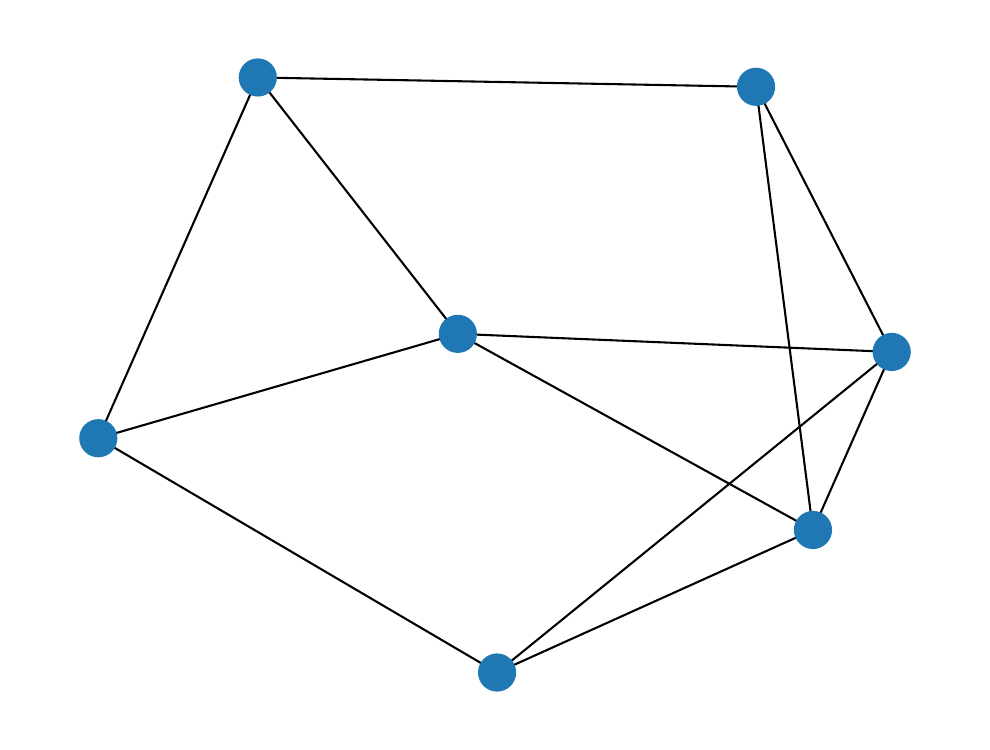} \ \ \ \ \ \
    \includegraphics[width=.3\textwidth]{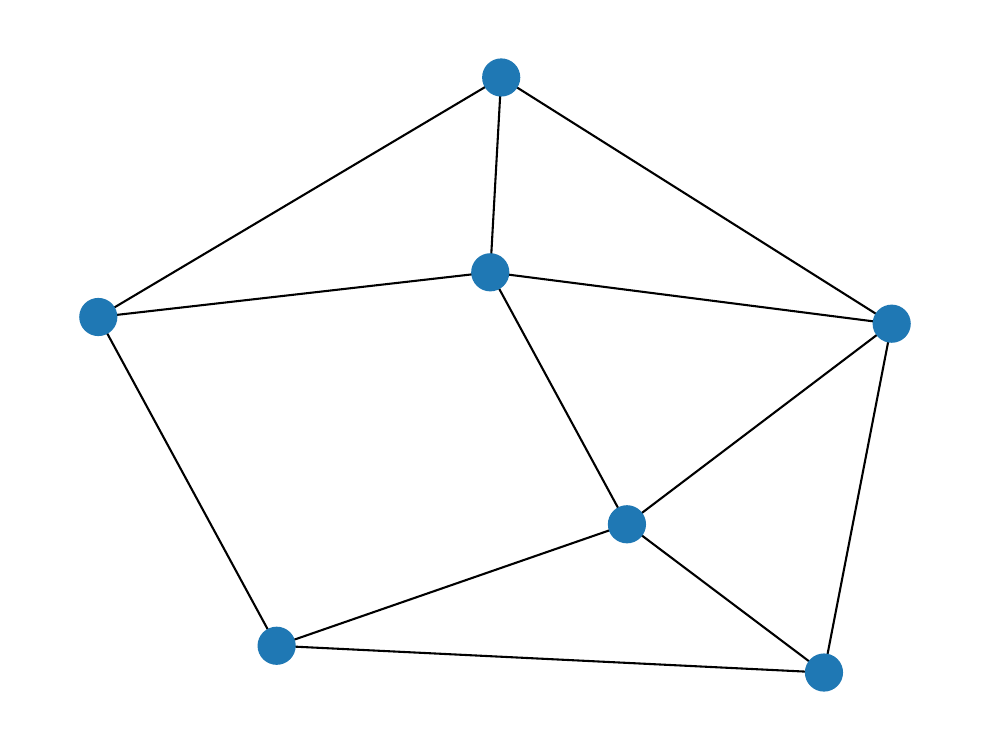}   
     \caption{One pair of graphs with the maximum similarity $s_2$ in family $N_7$.}
     \label{max_graphs7}
 \end{figure}

\section{Application in network theory}
\label{sec4}

In this section, we test our methodology on three different random network models in order to see weather it can discriminate among them. In particular, we use  the Erd\"{o}s-R\'{e}nyi model, the Barabási-Albert model (generating scale-free networks), and the Watts-Strogatz model (with ``small-world'' properties). We selected these three models since they are among the most commonly used models for generating random networks \cite{newman}.

We generate three graphs, $ER_1$, $ER_2$, and $ER_3$, in the Erd\"{o}s-R\'{e}nyi model using the parameters $n=100$ (the number of nodes in the graph) and $p=0.1$ (the probability that an edge exists between any pair of nodes). Next, graphs $BA_1$, $BA_2$, and $BA_3$ are generated by  the Barabási-Albert model  using the parameters $n=100$ (the number of nodes in the graph) and $m=3$ (the number of edges each new node creates when it is added). Finally,  graphs $WS_1$, $WS_2$, and $WS_3$ are generated by  the Watts-Strogatz model  using the parameters $n=100$ (the number of nodes in the graph), $k=6$ (each node is adjacent to 6 nearest neighbors in a ring topology), and $p=0.7$ (the rewiring probability). One example graph for each model is depicted in Figure \ref{random_net}. 

\begin{figure}[h!]
\centering
     \includegraphics[width=.4\textwidth]{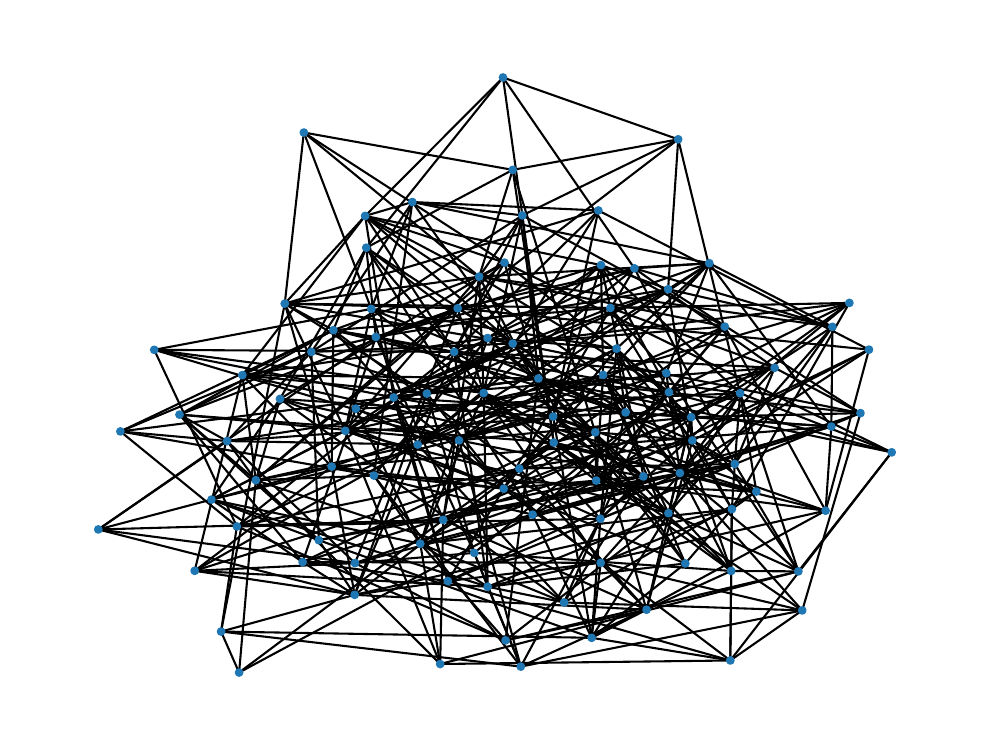} 
    \includegraphics[width=.4\textwidth]{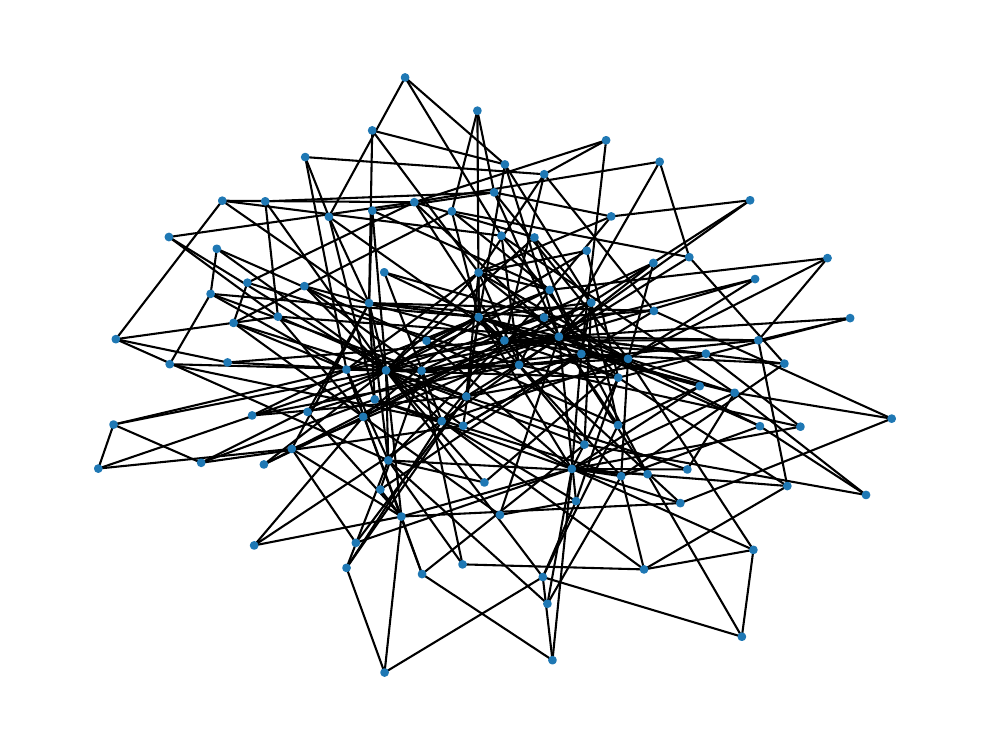}  
    \includegraphics[width=.4\textwidth]{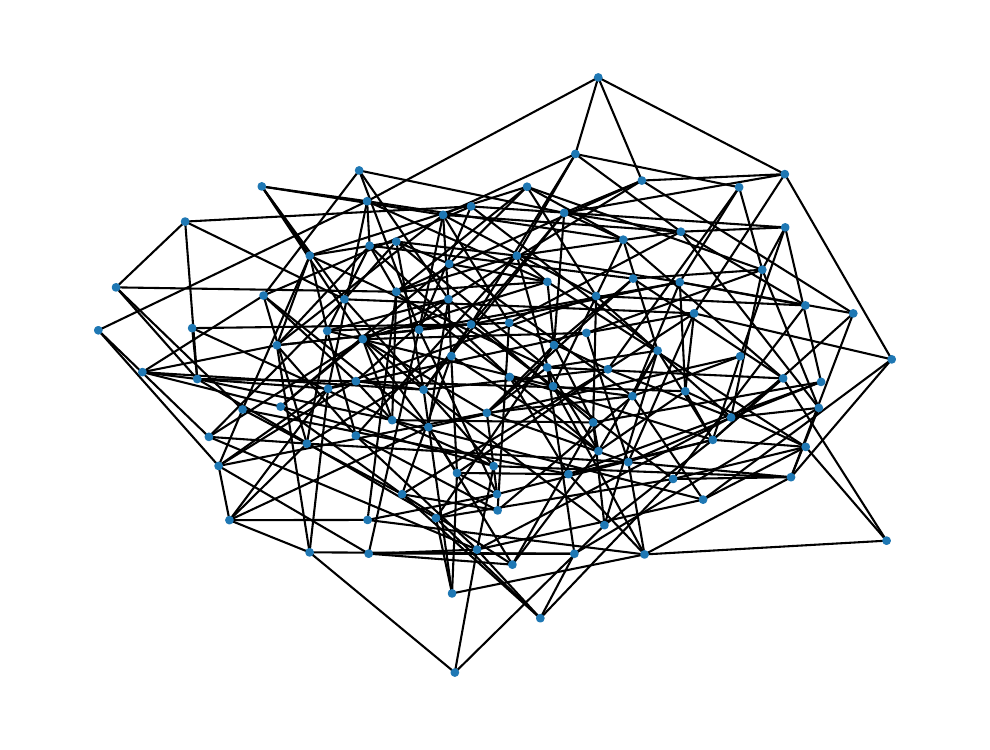} 
     \caption{Example networks generated by Erd\"{o}s-R\'{e}nyi model (upper left),  Barabási-Albert model (upper right), and Watts-Strogatz model (bottom).}
     \label{random_net}
 \end{figure}

 Then, we calculated the similarity $s_2$ for each pair of the considered random networks. Let $\mathcal{RN}$ be the set of all 9 considered networks. Moreover, if 
 \begin{eqnarray*}
     m &=& \min \{ s_2(G_1,G_2) \ | \ G_1,G_2 \in \mathcal{RN},\ G_1 \neq G_2 \} \ \textrm{ and}\\
     M &=& \max \{ s_2(G_1,G_2) \ | \ G_1,G_2 \in \mathcal{RN},\ G_1 \neq G_2 \},
 \end{eqnarray*}
 then for any pair of networks $G_1,G_2$ from $\mathcal{RN}$ we calculate the scaled similarity
 $$s_2' (G_1,G_2) = \frac{ s_2(G_1,G_2) - m}{M-m} $$
 to highlight the differences in the similarities. The results are gathered in Table \ref{random_scaled}.

\begin{table}[h!]
\begin{center}
\caption{\label{random_scaled} Scaled similarity $s_2'$ between random networks obtained by three different models.\\}
	\begin{tabular} {ccccccccc}\hline
    	&	$ER_1$	&	$ER_2$	&	$ER_3$	&	$BA_1$	&	$BA_2$	&	$BA_3$	&	$WS_1$	&	$WS_2$	\\ \hline \hline
$ER_2$	&	0.80	&		&		&		&		&		&		&		\\
$ER_3$	&	0.95	&	0.75	&		&		&		&		&		&		\\
$BA_1$	&	0.11	&	0.00	&	0.13	&		&		&		&		&		\\
$BA_2$	&	0.11	&	0.00	&	0.14	&	0.96	&		&		&		&		\\
$BA_3$	&	0.11	&	0.00	&	0.13	&	0.97	&	1.00	&		&		&		\\
$WS_1$	&	0.27	&	0.13	&	0.30	&	0.76	&	0.79	&	0.78	&		&		\\
$WS_2$	&	0.26	&	0.12	&	0.29	&	0.77	&	0.80	&	0.79	&	0.99	&		\\
$WS_3$	&	0.26	&	0.12	&	0.29	&	0.76	&	0.80	&	0.79	&	1.00	&	1.00	\\
	\hline
	\end{tabular}
    \end{center}
\end{table}

As one may see, our approach is able to distinguish among networks obtained by different models, namely the highest similarities are achieved for networks within the same model. On the other hand, significantly lower similarities are obtained for pairs of networks belonging to different models. For example, scaled similarities between Erd\"{o}s-R\'{e}nyi networks and Barabási-Albert networks are at most $0.14$. Slightly higher scaled similarities are obtained between Erd\"{o}s-R\'{e}nyi networks and Watts-Strogatz networks. On the other hand, scaled similarities between Barabási-Albert networks and Watts-Strogatz networks lie between $0.76$ and $0.80$. These findings show that our methodology may be applicable in network science for classification and comparison of networks. Namely, with further tuning and refinement, the approach can be adapted to capture the structural nuances of more complex networks, making it applicable to a broader range of real-world networks. 

\section{Application in chemistry}
\label{sec5}

In this section, we present results demonstrating the potential application of our method in chemistry, specifically in the assessment of molecular similarity. To evaluate its performance, we carried out a comparative analysis on a set of 75 decane isomers, comparing our approach with two established descriptors commonly used for molecular similarity evaluation. 
Specifically, the molecular structures of decane isomers were encoded using two types of molecular descriptors: Morgan circular fingerprints (with a radius of 2 and 2048 bits) \cite{rogers} and MACCS keys \cite{durant}. Then, the similarity between two molecules was calculated using the Tanimoto index \cite{ bajusz}, a widely accepted standard in molecular similarity computations. The correlations between the obtained similarities and the $s_2$ similarity for all pairs of decane isomers are depicted in  Figure \ref{chem_app}.

\begin{figure}[h!]
\centering
     \includegraphics[width=.5\textwidth]{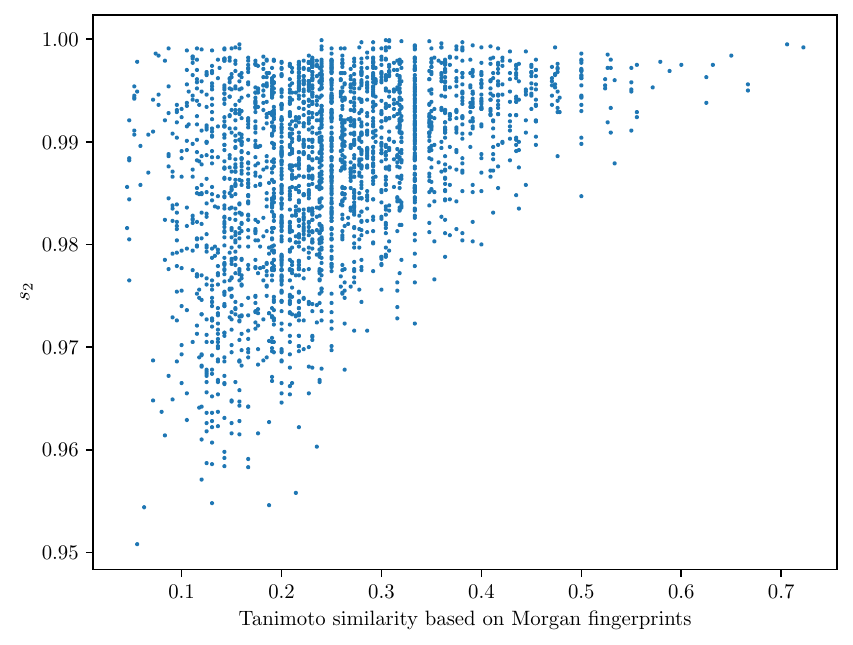} \hspace{-0.9em}
    \includegraphics[width=.5\textwidth]{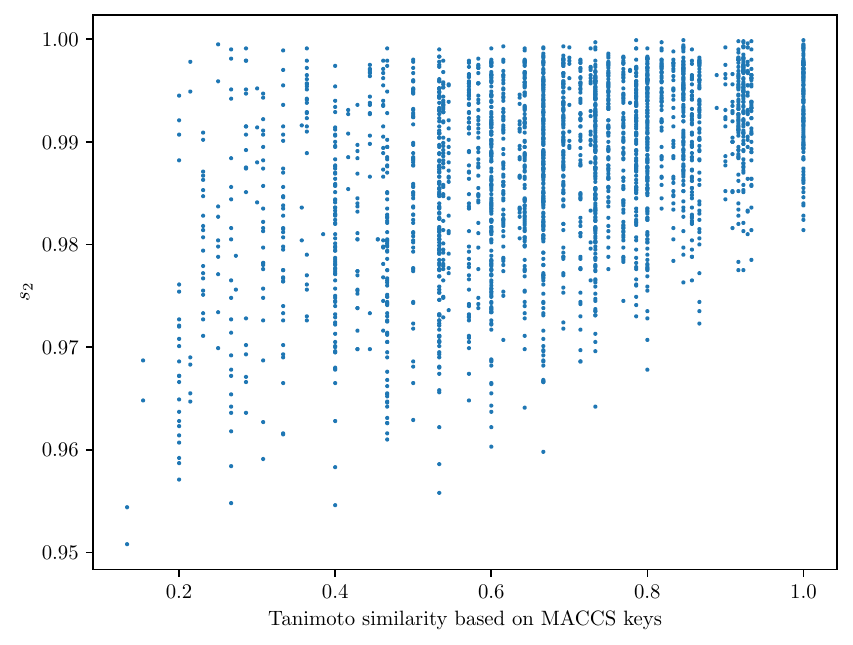}  
    \includegraphics[width=.5\textwidth]{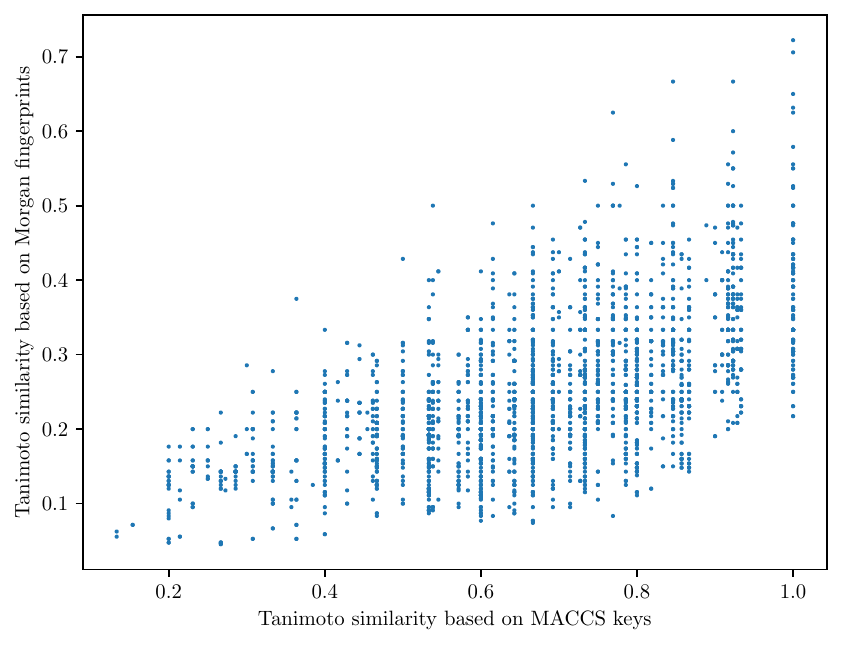} 
     \caption{Correlation between Tanimoto similarity based on Morgan fingerprints and $s_2$ (upper left), Tanimoto similarity based on MACCS keys and $s_2$ (upper right), and Tanimoto similarity based on MACCS keys and Morgan fingerprints (bottom), all for decane isomers.}
     \label{chem_app}
 \end{figure}

As observed from the results, the Tanimoto similarity based on Morgan fingerprints generally yields lower similarity scores for decane isomers compared to both MACCS key-based Tanimoto similarity and the $s_2$ structural similarity measure. This discrepancy likely arises from the nature of Morgan fingerprints, which emphasize local atomic environments and specific topological features. In contrast, MACCS keys capture broader, predefined substructural patterns that are more consistent across isomeric forms. Consequently, MACCS-based and $s_2$ similarities tend to align more closely.
It is also noticeable that both Tanimoto similarity measures, based on Morgan fingerprints and MACCS keys, tend to produce a relatively high number of molecule pairs with identical similarity scores. In contrast, such behavior is less apparent with the $s_2$ similarity values. This suggests that our method may offer increased sensitivity to subtle structural differences among isomers.
These findings indicate that our methodology may be well-suited for molecular similarity assessment. Furthermore, with the application of an appropriate weighting scheme to molecular graphs, the approach could potentially be extended to datasets containing more structurally diverse molecules.

 \section{ {Summary and conclusion}  }

In this paper, we introduced new graph similarity and distance measures. 
  While traditional methods such as graph edit distance and single-index topological measures provide valuable insights, they often face limitations in sensitivity or computational efficiency. The approach presented in this paper addresses these challenges by leveraging multiple topological indices to capture subtle structural differences while maintaining low computational cost. 
  
  Our results indicate that this method not only aligns well with established metrics but also demonstrates strong potential for large-scale applications and domain-specific analyses, including network modeling and molecular similarity studies. Specifically, the developed approach effectively distinguishes between real-world network models and measures molecular similarity with performance that exceeds that of established molecular similarity metric. Furthermore, applying suitable weighting schemes to both molecular graphs and network structures could be a promising way to extend the approach to datasets with more structurally diverse molecules and complex real-world networks. 

\section*{Funding information} 

\noindent Niko Tratnik and Petra \v Zigert Pleter\v sek acknowledge the financial support from the Slovenian Research and Innovation Agency: research programme No.\ P1-0297 and research  projects  N1-0285 (Niko Tratnik),  L7-4494 (Petra \v Zigert Pleter\v sek). Izudin Red\v{z}epovi\'{c} gratefully acknowledges financial support from the Serbian Ministry of Science, Technological Development, and Innovation (Agreement No. 451-03-137/2025-03/200252), as well as from the State University of Novi Pazar. Additional support was provided by the Faculty of Natural Sciences and Mathematics, University of Maribor, through the University of Maribor Visiting Researcher Programme 2025.  

\section*{Conflict of Interest Statement}

Not Applicable. The author declares that there is no conflict of interest.

\section*{Data availability statement}

All the computational results supporting the findings of this study are available within the paper.



\end{document}